\documentclass[a4paper,12pt]{article}
\usepackage[T2A]{fontenc}
\usepackage[koi8-r]{inputenc}
\usepackage{fullpage}
\usepackage{amssymb,amsmath}
\begin{document}

\title{Algorithmic information theory and martingales}

\author{Laurent Bienvenu, Alexander Shen}

\maketitle

\section{Introduction}

What is probability? What is (or should be) the subject of
probability theory? How this mathematical theory is (or should
be) applied to the ``real world''?

These questions were debated for centuries, and these
discussions go far beyond the scope of our paper. However, there
is a clear dividing line between two kinds of different
approaches; some of them attempt to define mathematically the
notion of an ``individual random object'' while the others move
this notion completely to the grey zone between ``pure''
probability theory (understood as a part of mathematics) and its
practical applications.

In practice, almost all mathematicians (and most
non-mathematicians), looking at the winning numbers of a lottery
for the last year and suddenly noticing that they are all even,
will conclude that something wrong happens. The same feeling
would arise if (as in the ``Rosenkrantz and Guildenstern are
dead'', the play by Tom Stoppard) the long sequence of heads
appears while tossing a (presumably fair) coin. However,
classical probability theory assigns to this sequence (say,
$100$ heads) the same probability $2^{-100}$ as to any other
sequence and does not try to explain why this sequence looks
``non-random'' and raises the suspicion.

This paradox (sequences with various regularities or symmetries
in them appear less random to us, even when each of them is just
as probable as any other outcome), occupied probabilists already in
the nineteenth century, including Laplace.\footnote{%
 ``C'est ici le lieu de d\'efinir le mot \emph{extraordinaire}.
 Nous rangeons, par la pens\'ee, tous les \'ev\'enements
 possibles en diverses classes, et nous regardons comme
 \emph{extraordinaires} ceux des classes qui en comprenement un
 tr\`es petit nombre. Ainsi, a jou de \emph{croix ou pile},
 l'arriv\'ee de \emph{croix} cent fois de suite nous parait
 extraordinaire, parce ques le nombre presque infini des
 combinaisons quit peuvent arriver en cent coups, \'etant
 partag\'e en s\'eries r\'eguli\`eres ou dans lesqulles nous
 voyone r\'egner un ordre facile \`a saisir, et en s\'eries
 irr\'eguli\`eres, celles-ci sont incomparablement plus
 nombreuses. La sortie d'une boule blanche d'une urne qui, sur
 un million de boules, n'en contient qu'une seule de cette
 couleur, les autres \'etant noires, nous parait encore
 extraordinaire, parce que nous ne formons que deux classes
 d'\'ev\'enement ordinaire, relatives aux deux couleurs. Mais la
 sortir du n$^\circ~475813$, par exemple, d'une urne qui
 renferme un million de num\'eros nous semble un \'ev\'enement
 ordinaire, parce que, comparant individuallement les num\'eros
 les uns aux autres, sans les partager en classes, nous n'avons
 aucune raisone de croire que l'un d'eux sortira plut\`ot que
 les autres.'' (``Essai philosophique sur les
 Probabilit\'es''~\cite{laplace}, VI~Principe).
 Peter G\'acs, who used this passage as an opening quote for his
 Dissertation~\cite{gacs-thesis}, comments: ``Laplace makes two
 informal suggestions (withouth strictly distinguishing them).
 First, he considers various classes of events, and views as
 extraordinary the small ones. (To make this precise, one would
 need to restrict attention to ``simple'' classes.) Second, he
 makes the assertion (without proof or even exact statement)
 that all outcomes of a given length having some regularity in
 them, grouped together, would still form a small class. (To
 make this precise, regularity must be defined
 appropriately.)''~\cite{gacs-private}}

However, the attempts to define \emph{mathematical} notions that
somehow capture the intuition of an individual random object (in
some idealized way) are not that old. Richard von Mises
suggestion (at the beginning of XXth century) was to base
probability theory on the notion of the so-called ``Kollektiv''
(an individual random sequence). These ideas were developed,
critically analyzed and made rigorous in 1930s by Wald, Ville
and Church (the latter gave a first precise definition of a
``random sequence'').

In 1960s and 1970s these notions were related to the notion of
complexity (amount of information, defined in algorithmic
terms), and now different definitions of randomness are well
studied in the framework of recursion theory and algorithmic
information theory.

In this paper we try to describe the main stages of this
development and its main achievements from the mathematical viewpoint
focusing on the role played by martingales.

This paper is based on published sources, discussion at the
Dagstuhl meeting (Seminar 06051, 29 January -- 3 February 2006;
C.~Calude, C.P.~Schnorr, P.~Vitanyi gave talks that were
recorded and made available at
\texttt{http://www.hutter1.net/dagstuhl} by Marcus Hutter) and
contributions of Leonid Bassalygo, Cristian Calude, Peter
G\'acs, Leonid Levin, Vladimir A. Uspensky, Vladimir Vovk,
Vladimir~Vyugin and others. It was initiated by Glenn Shafer
whose historical comments about Kolmogorov and Ville became a
starting point. (Of course, the people mentioned
are not responsible in any way for the
authors' flaws.)

\section{Collectives}

The first well known attempt to define mathematically the notion
of an individual random object was done by Richard von Mises in
his 1919 paper~\cite{mises-1919}. Then he elaborated his ideas
in the book published in 1928~\cite{mises-wsw}. He also made
some clarifying comments is his address delivered on September
11, 1940 at the meeting of the Institute of Mathematical
Statistics in Hanover, N.H. (USA) and published in
1941~\cite{mises-1941,mises-1941a}.

%Let us try to describe von Mises' ideas. First of all, one
%should have in mind that the modern standards of mathematical
%rigor were not so popular at that time, though gradually Mises
%had to clarify and formalize his ideas (probably under some
%pressure of others; one can feel this pressure
%reading~\cite{mises-1941,mises-1941a}).

Mises explains that probability theory studies a special class
of natural phenomena, like tossing a coin, rolling a dice, or
other repetitive experiments. Geometry tries to capture and
axiomatize the real-world notion of space; in a similar way
probability theory captures and axiomatizes the properties of
random phenomena, called ``collectives'' (German: Kollektiv)
in Mises' paper.
Informally speaking, collectives are (according to Mises)
plausible sequences of outcomes we can get by performing
infinitely many independent trials of some experiment. He
formulated two axioms for the notion of collectives. For
simplicity, we state them for a collective with two values, e.g.,
the sequence of heads and tails obtained by coin tossing (where
the coin is potentially unbalanced, i.e., the outcome ``tails''
may appear more (or less) often than ``heads''):

\medskip

\textbf{I}. There exists a limit frequency: if $s_N$ is the
number of heads among the first $N$ coin tosses, the ratio
$s_N/N$ converges to some real $p$ as $N\to\infty$.

\smallskip

\textbf{II}. This limit frequency is stable: if we select a
subsequence according to some ``selection rule'', then the
resulting subsequence (if infinite) has the same limit frequency.

\medskip

Axiom I is quite natural: if we want to explain informally what
probability is, we say something like ``repeat the experiment
many times until the frequency of some event (say, head on a
coin) becomes almost stable; this stable value is called a
probability of the event''.

What is the second axiom needed for? Remember that collectives
should represent \emph{plausible} sequences of outcomes of
independent trials. Suppose somebody tells you that
flipping a coin produced the sequence
        $$
0101010101010101010101010101\ldots
        $$
where $0$ (heads) and $1$ (tails) alternate. Would you believe
this? Probably not. Globally, the limit frequency of $0$ and $1$
in this sequence exists and is equal to $1/2$. But this
sequence does \emph{not} look plausible as a sequence of outcomes,
as it presents some
highly suspicious regularity. This is where axiom II comes into
place: if one selects from this sequence the bits in even
positions, one gets a new sequence
        $$
1111111111111111111111111111\ldots
        $$
in which the frequency of ones is different ($1$ instead of
$1/2$).

Probability theory, according to Mises, needs to define its
subject, and this subject is the properties
of collectives and operations that transform collectives into
other collectives. Mises uses the following example: take a
collective (a sequence of zeros of ones) and cut it into 3-bit
groups. Then replace each group by an individual bit according
to the majority rule. Probability theory has to find the limit
frequency of the resulting sequence if the limit frequency of
the original one is known.

In his early papers Mises explained in quite informal way which
selection rules are allowed: the selection rule should decide
whether a term is selected or not, using only the values of the
preceding terms but not the value of the term in question. For
example, selection rule may select terms whose numbers are
prime, or terms that immediately follow heads in the sequence,
but not the terms that \emph{are} heads themselves.

The existence of collectives, according to von Mises, is an
observation confirmed by our experience, e.g., by thousands of
people who invented different systems to beat the casino but all
failed in the long run (principle of ``ausgeschlossenen
Spielsystem'', as Mises said).

\section{Clarifications. Wald's theorem}

Of course, Mises' approach was quite vulnerable from the
mathematical viewpoint. What is a selection rule? Do collectives
exist at all?

Answering these objections, Mises adopted a more formal
definition of a selection rule suggested by A.~Wald
(see, e.g.,~\cite{wald} and~\cite{mises-1941}).
Assume for simplicity that a sequence is formed by zeros and
ones. The selection rule is a total function
$s\colon\{0,1\}^*\to\{0,1\}$. Here $\{0,1\}^*$ is a set of all
finite binary strings. Applying selection rule $s$ to an
infinite binary sequence $\omega_1\omega_2\ldots$ means that we
select all terms $\omega_i$ such that
$s(\omega_1\omega_2\ldots\omega_{i-1})=1$; the selected terms
are listed in the same order as in the initial sequence.

The condition II for a selection rule $s$ says that for a
collective the selected subsequence either should be finite or
should have the same limit frequency as the entire sequence.

Therefore we get a formal definition of a collective as soon we
fix some class of selection rules. The evident problem here is
that if we consider \emph{all} selection rules of the described
type, collectives (non-trivial ones, with limit frequency not
equal to $0$ or $1$) do not exist. Indeed. for every set $S$ of
natural numbers there exists a selection rule that selects the
terms $\omega_i$ for $i\in S$ (the function $s$ depends only on
the length of its argument). Using for a given sequence
$\omega_1\omega_2\ldots$ the set $S$ of all $i$ such that
$\omega_i=0$ (or $\omega_i=1$), we get a contradiction.

Wald~\cite{wald}\footnote{%
      A short note without proofs was published
      earlier~\cite{wald-1936}.}
provided a kind of solution
for this problem.
He proved that for any \emph{countable} family of selection
rules and for any $p\in(0,1)$ there is a continuum of sequences
that satisfy the axioms I (with limit frequency~$p$) and II for
this class of selection rules.

Today this statement looks almost trivial: indeed, if a given
selection rule $s$ is applied to a $B_p$-randomly chosen
sequence, where $B_p$ is Bernoulli distribution with
parameter~$p$, the selected subsequence has the same
distribution $B_p$, so the Strong Law of Large Numbers
guarantees that the set of sequences that do not satisfy II for
a given $s$ has $B_p$-measure zero; the countable union of null
sets is a null set and its complement has continuum cardinality.

However, Wald wanted to give a constructive proof of this
result; Theorem V (\cite{wald}, p.~49) says that if a
``konstruktiv definiertes abz\"ahlebare System von
Auswahlvorschriften'' is given, ``so kann man Kollektiv
$\langle\ldots\rangle$ konstruktiv definieren'' (if a countable
system of selection rules is defined constructively, there
exists a constructively defined collective).

Note that there is
no formal definition of ``constructive'' objects in Wald's paper;
he just provides a construction of a collective that refers
to selection rules (uses them as an oracle, in modern terminology).
The collective sequence is constructed inductively. Let us
explain the idea of the construction in a simple case when only
finitely many selection rules $s_1,\ldots,s_n$ are considered
and sequence of zeros and ones has limit frequency $1/2$.

At the $i$th step of the construction we should decide whether
$\omega_i$ is $0$ or $1$. At that time we already know which of
the rules $s_1,\ldots,s_n$ would include $\omega_i$ in the
selected subsequence. In other terms, we know a Boolean vector
of length $n$. The entire sequence (that we have to construct)
would be therefore split into $2^n$ subsequences that correspond
to $2^n$ values of this Boolean vector. Now the main idea: each
of these $2^n$ sequences should be $0101010\ldots$ (zeros and
ones alternate starting with zero). This determines the sequence
$\omega$ uniquely. Since $\omega$ is a mixture of $2^n$
sequences that have limit frequency $1/2$, the entire sequence
$\omega$ has the same limit frequency.

And if we apply selection rule $s_i$ to $\omega$, we get a
mixture of $2^{n-1}$ of these subsequences (corresponding to
$2^{n-1}$ Boolean vectors where $s_i$ is playing). Each sequence
has limit frequency $1/2$, and their mixture has therefore the
same limit frequency.

In fact the construction for countably many selection rules is
quite similar: we just have to add new rules one by one when the
sequence is so long that the boundary effects cannot destroy the
limit frequency.

In fact Wald proves more: he considers not only the two-element
set $\{0,1\}$, but any finite set (Theorem I, p.~45). Then he
considers the case of infinite set $M$ (Theorem II--IV,
pp.~45--47; we do not go into details here, but to get a
reasonable definition of a collective for infinite $M$ one
should either consider countable $M$ or a restricted class of
events). Theorems V--VI (p.~49) observe that the resulting
collectives are ``constructive''.

Based on Wald's results, Mises~\cite{mises-1941} concludes that
the notion of colletive can be studied without contractictions:
we can consider all the selection rules we want to use and their
combinations; though we do not know them in advance, one may
reasonably assume that they form a
finite or
countable set and therefore collectives (with respect to this
set) do exist.

Wald's results show, in a sense, that the requirements I and II
are not too strong. But other objections to the notion of
collective, raised by Ville in his book~\cite{ville-etude}, say
that these requirements are too weak: not only collectives
exist, but one can construct some collective in the sense of
Mises' definition that does not look random.

\section{Ville's objections. Martingales}
\label{martingale}

Let us explain Ville's objections. The requirement II can be
reformulated in terms of games as follows. (For simplicity we
consider the case when limit frequency is $1/2$.) A player comes
into a casino where a coin is tossed infinitely many times, and
can (for each tossing) decide to make a bet or to skip it
depending on the results of a previous tossings (according to
the selection rule she has in mind). Her initial credit is \$0,
and she is allowed to incur arbitrarily large debts. All bets
are for the same amount of money, say \$1, which the player
loses or doubles, depending on whether her guess was correct or
not. Let $c_N$ be the player's capital after $N$ games. The
player wins (after infinitely many games) if she makes
infinitely many bets and the ratio $c_N/N$ does not converge to
zero.

(This game deviates from the original idea of a selection rule:
instead of just choosing of a subsequence we are allowed also to
reverse some of the terms chosen. However, this gives an equivalent
definition since we may consider separately the ``positively'' and
``negatively'' chosen terms; if both subsequences have limit
frequencies $1/2$, the ratio $c_N/N$ does converge to
$0$. Note also that this definition assumes that the coin is fair.)

We have reformulated Mises' definition in terms of a game,
but this game looks rather
unnatural. Yes, for a ``really random coin'' we would expect
that $c_N/N$ converges to $0$ (at least after we learned the
strong law of large numbers). But is it the only thing we would
expect? Imagine, for example, that $c_N$ is always positive and
goes slowly but steadily to infinity, so $c_N/N\to 0$ but
$c_N\to+\infty$. This would mean that the player manages to make
arbitrarily large amounts of money without incurring debts. In
that case, would we agree with the assumption that she is
playing with a fair coin?

Ville suggested a different kind of gambling games,
which are much more natural. In his games we come to the casino
with some fixed amount of money (say, \$1) and can use it (in
whole or in part) for betting, but cannot go negative. In other
terms, if we have $s$ before the next game, we can bet any
amount $s'\le s$ on zero or one. If our guess is incorrect, the
money is lost, and our capital becomes $s-s'$, otherwise the
money is doubled, and our capital is then $s+s'$.

Mathematically such a strategy is represented by a function $m$
whose arguments are finite binary strings and values are
non-negative reals. The value $m(\omega_1\ldots\omega_n)$ is our
capital after we have played $n$ times getting outcomes
$\omega_1,\ldots,\omega_n$; the value $m(\Lambda)$ (where
$\Lambda$ denotes the string of length zero) is the initial
capital, which we assume to be positive. The rules of the game
dictate that
        $$
m(x)=\frac{m(x0)+m(x1)}{2}   \eqno(*)
        $$
Here $x$ is some binary string (representing some moment in the
game), $x0$ and $x1$ are obtained by adding $0$ or $1$ to $x$
and correspond to two possible outcomes in the next round. The
requirement says that $m(x)$ is the average between two
possibilities, i.e., our possible gain and loss are balanced.
Ville used the name \emph{martingale} for functions that have
property $(*)$. (One may also allow the martingales to have negative
values, but we use only non-negative martingales in the sequel.)

A martingale $m$ (i.e., the player that uses
corresponding strategy) \emph{wins} against a sequence
$\omega_1\omega_2\ldots$ if the values
$m(\omega_1\omega_2\ldots\omega_n)$ are unbounded. Now we can
switch from Mises' selection rules to martingales and say that a
sequence $\omega=\omega_1\omega_2\ldots$ is a collective (in a
new sense) if all martingales from some (countable) family do
not win against~$\omega$.

\medskip

To support this change in the class of games, Ville notes that:

\begin{itemize}
        \item
Martingales provide a generalization of Mises' games (with limit
frequency $1/2$): for any selection rule one can construct a
martingale that wins against every sequence that does not
satisfy axiom II when this selection rule is applied.
        \item
The notion of martingale matches well the notion of a null set
(set of measure~$0$) used in classical probability theory:
for every martingale $m$, the set of all sequences
against which~$m$ wins is a null set (has measure $0$) according
to the uniform Bernoulli distribution.
        \item
The reverse statement is also true: for every null subset
$X\subset\{0,1\}^\infty$ there exists a martingale $m$ that wins
against every element of $X$. (Together with the strong law of
large number this implies the first statement in the list).
\end{itemize}

(The proofs are quite natural: first we prove the finite
versions of these results saying that (\textbf1)~the probability
to transform initial capital $1$ into some $C$ during $N$ games
does not exceed $1/C$; (\textbf2)~for every $N$ and for every
set of $N$-bit sequences that contains $\varepsilon$-fraction of
all sequences of length~$N$, there is a martingale that wins
$1/\varepsilon$ on every sequence from this set.)

\medskip

Martingales have some other nice properties. One may ask why our
winning condition says that martingale is unbounded: isn't it
more natural to require that its values tend to $+\infty$ (a
strong winning condition)? The answer is that it
does not matter much, as the
following simple observation shows: for every martingale $m$
there exist another martingale $m'$ that strongly wins against
a sequence $\omega$ if $m$ wins against~$\omega$. (The
martingale $m'$ should save part of the capital when the capital
reaches some bound and use only the remaining part for playing,
waiting until it has enough to save again, etc.)

Another nice property is the possibility of combining
martingales: if $m_i$ are arbitrary martingales, the weighted
sum $\sum_i \alpha_i m_i$ (where $\alpha_i$ are some positive
reals with sum~$1$) is a martingale that wins against a sequence
$\omega$ if and only if at least one of $m_i$ wins
against~$\omega$. (Recall that we consider only non-negative
martingales.)

\section{Ville's example}

The arguments above may convince you that martingales have more
nice properties than just selection rules.\footnote{In fact,
   at Ville's time these arguments did not sound very convincing
   even to some experts: W.~Feller wrote in his Zentralblatt
   review of one of the first Ville's papers: ``Aus
   unerfindlichen Gr\"unden will nun Verf. den Auswahlbegriff so
   ab\"andern (``martingale'' statts Auswahl) da\ss\ jede
   Nullmenge als Ausnamemenge bei passendem $S$ autreten kann'',
   both reproducing the main argument of Ville (the possibility
   to exclude any null set) and finding it unconvincing
   (``unerfindlichen Gr\"unden''), see~\cite{shafer}.}
But is this difference
essential? If we switch from selection rules to martingales,
do we get stronger requirements for
random sequences (collectives)? Ville showed that it is indeed
the case, proving the following result.

\begin{quote}
For any countable family $\mathcal{S}$ of selection rules there
exists a sequence $\omega$ that satisfies requirement II (with
limit $1/2$) when rules from $\mathcal{S}$ are used but every
prefix of $\omega$ has at least as many zeros as ones
(\cite{ville-etude}, p.~63, Remarque).
\end{quote}

(In fact, Ville proved more;
Theorem~4, p.~55, provides also some bounds for the speed of
convergency.)

This proof raises a historical question. In fact, Ville's
argument is very close to Wald's argument used in~\cite{wald}:
the sequence is splitted into subsequences and inductive
construction is performed; Wald does not discuss the one-sided
convergence explicitly, but it is obtained in a straightforward
way as a byproduct of Wald's conctruction. Indeed, let us say
that a sequence is ``biased'' if every prefix has at least as
many zeros as ones (frequency of ones does not exceed $1/2$). If
we merge biased sequences, the result is also a biased sequence;
note also that the sequence $01010101\ldots$ is biased.

However, Ville does not mention this similarity (though Wald's
paper is mentioned many times in Ville's book and the existence
result is quoted with reference to Wald). It is especially
strange since the explanations given in Wald's paper are quite
clear~--- probably more clear than Ville's argument, which is
written in a rather technical way. May be this heavy technical
style of Ville's paper was the reason why other authors prefer
to give their own reconstruction of the proof instead of
following the details of Ville's paper (see, e.g.,~\cite{lieb}
and references within).

\section{More about Ville's example}

Establishing the difference between selection-based and martingale-based
definitions of randomess,
Ville also showed that there is a martingale that wins against
every ``biased''
sequence (a sequence whose
prefixes have more zeros than ones).
This is a consequence of the law of
iterated logarithm; it implies that the set of all biased sequences
has measure zero, so we can use the results mentioned in
Section~\ref{martingale}. However, let us provide a simple direct
construction of such a martingale just for illustration.

Let $\omega$ be a binary sequence; let $d_n$ be the difference
between the numbers of zeros and ones in $n$-bit prefix of
$\omega$. We assume that the difference $d_n$ is always
non-negative. The limit $d=\liminf d_n$ is then also
non-negative; it can be finite or $+\infty$.

It is easy to construct a martingale that wins against any
biased sequence with ${d=+\infty}$. Imagine that you come into a
casino knowing in advance that (1)~the number of ones never
exceeds the number of zeros and (2)~the difference between them
tends to infinity. How can you become infinitely rich? Just bet
a fixed amount (not exceeding the initial capital) at every
step. The condition (1) guarantees that you will never go
negative and always have enough money to bet; the condition (2)
guarantees that your capital tends to infinity.

Now assume that the casino sequence is biased and $d$ is finite.
How can you win then? In this case the difference goes below $d$
only finitely many times, and starting from some time $T$ it is
at least $d$ being equal to $d$ infinitely many times. A
conclusion: if you see (after the initial period of length $T$)
that the difference is $d$, you know that the next coin tossing
provides a head, so you bet on it with no risk. This allows
you to become infinitely rich if you know $d$ and $T$ in
advance.

So we have one martingale $m$ that wins against any biased
sequence with $d=+\infty$ and a countable family $m_{d,T}$ of
martingales who win against sequences with given $d$ and~$T$.
As we have noted, this countable family of martingales can
be combined into one martingale.

\medskip

There is a large variety of possible interpretation of Ville's
example. One can treat this example as a failure of Mises'
approach: it shows that requirements I and II that guarantee
frequency stability (and therefore establish the very notion
of probability) are not strong enough to provide a
satisfactory definition of
a random sequence (collective): a martingale cannot win against a
``real coin'' but still can win against a collective formally
defined in terms of selection rules.

One may say also that axioms I and II do not pretend to capture
all properties of ``really random'' sequence but only some of
them needed to define the notion of probability, and therefore
the Mises' notion of collective can be considered as an upper
bound for the class of ``really random'' sequences.

Finally, one can say also that replacing selection rules by a
stronger martingale requirement, we harmonize the idea of a
random sequence with the measure-theoretic understanding of laws
of probability theory, therefore giving new life to Mises'
approach and getting a better notion of randomness.

It would be interesting to reconstruct the real attitude of
Mises, Ville, Frechet and others; however, this again goes far
beyond the scope of the article. Let us note nevertheless that
the only place where Ville is mentioned in~\cite{mises-1964} has
nothing to do with martingales (it is a paper on game theory).
Things become even more complicated when we try to interpret
Mises' remark in~\cite{mises-1919} when he says: ``Solange man
etwa nur die Zahlen $1$--$10000$ betrachtet, bietet die
Anordnung der Ziffern an der 5. Stelle [in the table of
logarithms] tats\"achlich das ungef\"ahre Bild eines empirisches
Kollektivs und man kann auch die S\"atze der
Wahrscheinlichkeitsrechnung n\"aherungsweise darauf anwenden.''
This quote shows that for him (at least at that moment) the
behavior of the 5th decimal digit in the table of logarithms of
integers $1$--$10000$ looks like ``empirical collective'' and
this sequence satisfies the laws of probability theory to a
certain extent (while for bigger numbers the regularities show
up). Note that logarithms are computable, so there exists a
computable selection rule that selects only zeros from this
sequence. One may speculate that Mises had in mind some notion
of ``pseudorandom'' sequence that satisfies the axiom II only
for simple enough selection rules, but this remark remains
isolated in his paper and it is hard to say what he really
meant.

\section{Church definition of randomness}

Approximately at the same time, in 1930s, a theory of
computable functions was developed by Kleene, Church, Turing and
others. It provided a very natural class of selection rules:
computable rules, where the function
$s\colon\{0,1\}^*\to\{0,1\}$ is a total computable function.
This class contains almost all rules we can think of; it also
has nice closure properties needed to prove theorem about
collectives. For example, it is closed under composition, and
this can be used to prove that a sequence obtained from a
collective by a selection rule is again a collective.

This step (combining recursion theory with Mises' approach) was
done in 1940 by Church~\cite{church-random}: he called a
sequence random if it has limiting frequency and, moreover, any
computable selection rule produces either finite sequence or a
sequence with the same limit frequency.

In fact, Church could do the same with Ville's definition and
define random sequences using computable martingales. But
probably he did not realize the importance of martingales.

More details about the evolution of the randomness notion from
Mises to Church can be found in a historical survey of
Martin-L\"of~\cite{martin-lof-1969}.

\section{An intermission}

In the 1940s and 1950s the notion of an individual random
sequence did not attract much attention. At that time the
measure-theoretic approach to probability theory became
gradually more and more popular (and, in particular, the notion
of martingale was embedded into the framework of measure
theory).

Another important change during these 20 years was the
development of the theory of computation. In 1930s theory of computation
appeared as a kind of exotic thing developed by logicians
that is
using
strange tools like recursive functions (with quite unnatural
definition), $\lambda$-calculus (even more peculiar definition)
or fictional devices called ``Turing machines''. But after
twenty years the notion of a computer program became quite
familiar; many mathematicians played with computers (i.e.,
programmed them~--- computer games for dummies were almost
unknown at that time) as a part of their job or just for fun.

This prepared a next step in the development of randomness
notion when the connections with the complexity
(incompressibility) was understood.

\section{Complexity and randomness in 1960s}

Recall the question we started with: why does the long sequence
of zeros (heads) look suspicious while the other sequence of the
same length (having the same probability $2^{-n}$ according to
the classical theory) looks OK? What is the difference between
these two sequences?

Now, when the notion of computer program became familiar, the
difference between them is evident: the first sequence (zeros)
can be generated by a short program while the other one
(non-suspicious) cannot.

So there is no surprise that the ideas of complexity of a finite
object (defined as the length of a shortest program that
generates this object) were developed independently in different
places and by different people. This kind of complexity is often
called \emph{description} complexity, as opposed to
\emph{computation} complexity, since we ignore the time needed
to generate an object and look only at the length of the
generating program.

There were other (not related to randomness) reasons to consider
description complexity. One of these reasons was the quantitative
analysis of undecidability. ``Undecidable algorithmic problems
were discovered in many fields, including algorithms theory,
mathematical logic, algebra, analysis, topology and mathematical
linguistics. Their essential property is their generality: we look
for an algorithm that can be applied to every object from some
infinite class and always gives a correct answer. This general
formulation makes the question not very practical. A practical
requirement is that algorithm works for every object from some
finite, though probably very large, class. On the other hand,
the algorithm itself should be practical. $\langle\ldots\rangle$
Algorithm is some instruction, and it is natural to require that
this instruction is not too long, since we need to invent this
algorithm\ldots\ So an algorithmic problem could be unsolvable
in some practical sense even if we restrict inputs to some
finite set'' (A.A.~Markov~\cite{markov-1967}, p.~161; this paper
provides proofs for the results announced in~\cite{markov-1964})

Note also that the idea of measuring the complexity of a message
as the length of its shortest ``encoding'' was quite familiar
due to Shannon information theory (though the encodings
considered there are very restricted).

Earlier (in~\cite{solomonoff-1964a,solomonoff-1964b}; these
papers are based on technical reports that go back to 1960 and
1962) R.~Solomonoff considered similar notions in the context of
inductive inference (somebody gives us a long sequence; we want
to know what is the reasonable way to predict the next term of
this sequence knowing the preceding terms).

G.~Chaitin~\cite{chaitin-home} tells that entering a Bronx
High School of Science (in 1962) he wrote an essay where the
idea of randomness as an absence of short description was
mentioned; later, in 1965, after his first year in City College,
he wrote a paper that was submitted to the Journal of the ACM
and finally published in two
parts~\cite{chaitin-1966,chaitin-1969}. In~\cite{chaitin-1966}
he defines a complexity measure of a binary string in terms of
the size of a Turing machine; in~\cite{chaitin-1969} the
complexity is defined in more general terms (in the same way as
in Kolmogorov paper~\cite{kolmogorov-1965}, see
below).\footnote{%
        The most famous discovery of Chaitin is probably the
        proof of G\"odel incompleteness theorem based on the
        Berry paradox~\cite{chaitin-1971}; we don't discuss it
        here.}

L.A.~Levin~\cite{levin-interview,levin-memoir} tells that being a student of
a high school for gifted children in Kiev (USSR, now Ukraine) in
1963/4, he was thinking about the length of the shortest
arithmetic predicate that is provable for a single value of its
parameter but did not know how make this definition invariant
(how to make the complexity independent of the specific
formalization of arithmetics). Next year (1964/1965) he moved to
Moscow where a special boarding school for gifted children was
founded by A.~Kolmogorov, and told about this idea to
A.~Sossinsky who was at that time a teacher in this school.
Sossinsky asked Kolmogorov and Kolmogorov replied that in one of
his forthcoming papers this question was answered.\footnote{%
        Here is the Russian quotation from~\cite{levin-memoir}:
        ``Тема, которой Андрей Николаевич тогда увлекался~---
        общие понятия сложности, случайности, информации~---
        волновала меня чрезвычайно. Как многие молодые люди,
        я искал самых фундаментальных концепций. Но такие
        ``первичные'' теории, как логика или теория алгоритмов,
        смущали меня своей ``качественной''  природой~---
        там нечего было ``посчитать''. На самом деле, я ещё
        в Киеве пытался дать определение сложности (я называл
        её ``неестественность''), но не мог доказать её
        инвариантности. В Москве я рассказал о своих неудачах
        Сосинскому, он спросил Колмогорова и принёс мне
        поразительный ответ: Колмогоров как раз доказал то,
        что я не смог и уже вот-вот выйдет его подробная
        статья! Тогда я решил во что бы то ни стало поступить
        в МГУ и стать учеником Андрея Николаевича.''}

This was the paper~\cite{kolmogorov-1965} that soon became the
main reference for the definition of complexity; now the
complexity defined as the length of the shortest program is
often called ``Kolmogorov complexity''. The paper was called
``Three approaches to the quantitative definition of
information'', and one of the approaches (the algorithmic one)
defined the complexity of a binary string as the length of the
shortest program producing it, assuming the programming language
is optimal, and proves the existence of such an optimal language
(for the technical details see the paper or any of the tutorials
on Kolmogorov complexity, e.g.,~\cite{shen-notes}).

This Kolmogorov paper had several historical reasons to become
most popular (among many expositions of the same ideas,
including the above mentioned). It was the first publication
where the rigorous definition of complexity was given and
universality theorem was proved. (This was done also in the
second part of Chaitin's article submitted in November 1965,
after Kolmogorov's publication, and published only in 1969.
Solomonoff's papers did not contain an explicit definition of
complexity.)

Second, Kolmogorov was famous as one of the greatest
mathematicians of his time, and therefore his papers attracted a
lot of attention. And being one of the founders of probability
theory, he has a clear vision of the role that complexity can
play in the foundations of probability theory (in the definition
of individual random object and in information theory). So his
paper was concise and well written.%
        \footnote{%
Chaitin's papers start with a lot of technical details related
to the counting of Turing machines states. Solomonoff's
paper~\cite{solomonoff-1964a} contains passages like ``The
author feels that Eq.~(1) is likely to be correct or almost
correct, but that the methods of working the problems of
Sections 4.1 to 4.3 are \emph{more} likely to be correct than
Eq.~(1). If Eq.~(1) is found to be meaningless, inconsistent or
somehow gives results that are intuitively unreasonable, then
Eq.~(1) should be modified in ways that do not destroy the
validity of the methods used in Sections 4.1 to 4.3''~--- not
very encouraging for the readers, to say the least. Levin
remembers that when he was instructed by Kolmogorov to read and cite the
work of Solomonoff, he was frustrated by this kind of attitude and
soon gave up.

Section 3.2.1 of~\cite{solomonoff-1964a} contains the following
sentence: ``Although a proof [of some statement, related to a
definition called Eq.~(1); this definition contained an error,
as Solomonoff found later] is not available, an outline of the
heuristic reasoning behind this statement will give clues as to
the meanings of the terms used and the degree of validity to be
expected of the statement itself''. But later in the same
paragraph a very clear proof of universality theorem is provided
for the readers who are not confused by previous remarks and are
able to extract its statement out of the proof. This paper also
contained a lot of other ideas that were developed much later;
e.g., in Section 3.2 Solomonoff gives a nice simple formula for
predictions in terms of the conditional a priori probability,
using monotonic machines much before Levin and Schnorr. (In 1978
Solomonoff formally proved that this formula works for all
computable probability distributions,
see~\cite{solomonoff-1978}.)

In fact, Solomonoff's main interest was inductive inference. He
tried to formalize the ``Occam's Razor'' principle in the
following way: base your prediction on the simplest ``law'' that
fits the data, say the simplest program that could generate it.
This requires a definition of ``simplecity'', and it was in this
context that Solomonoff defined complexity in terms
of description length
and proved its invariance. (His actual prediction
formula uses conditional a priori probability, based on all
possible programs that fit the data, with longer programs
entering with smaller weights.) }
Therefore it is no wonder that among many people who came to
very close ideas, Kolmogorov got the most attention.\footnote{%
   When Kolmogorov has came to the definition of complexity?
   In his 1963 paper~\cite{kolmogorov-1963} Kolmogorov makes
   some remarks that partially explain how he came to the
   complexity notion: ``I have already expressed the view
   $\langle\ldots\rangle$ that the basis for the applicability
   of the results of the mathematical theory of probability to
   real `random phenomena' must depend on some form of the
   frequency concept of probability, the unavoidable nature of
   which has been established by von Mises in a spirited manner.
   However, for a long time I had the following views:

   (1)~The frequency concept based on the notion of limiting
   frequency as the number of trials increases to infinity, does
   not contribute anything to substantiate the applicability of
   the results of probability theory to real practical problems
   where we have always to deal with a finite number of trials.

   (2)~The frequency concept applied to a large but finite
   number of trials does not admit a rigorous formal exposition
   within the framework of pure mathematics.

   Accordingly I have sometimes put forward the frequency
   concept which involves the conscious use of certain not
   rigorously formal ideas about `practical reliability',
   `approximate stability of the frequency in a long series of
   trials', without the precise definition of the series which
   are `sufficiently large'\ldots\

   I still maintain the first of the two theses mentioned above.
   As regards the second, however, I have come to realize that
   the concept of random distribution of a property in a large
   finite population can have a strict formal mathematical
   exposition. In fact, we can show that in sufficiently large
   populations the distribution of the property may be such that
   the frequency of its occurrence will be almost the same for
   all sufficiently large sub-populations, when the law of
   choosing these is sufficiently simple. Such a conception in
   its full development requires the introduction of a measure
   of the complexity of the algorithm. I propose to discuss this
   question in another article. In the present article, however,
   I shall use the fact that there cannot be a very large number
   of simple algorithms.'' In this quote Kolmogorov suggested a
   finitary Mises-style approach that uses selection rules of
   bounded complexity, but does not explain what complexity is;
   also he does not speak here about definition of randomness in
   terms of complexity (directly, without using selection
   rules).

   Asked when Kolmogorov came to his definition of complexity,
   Martin-L\"of writes~\cite{martin-lof-interview}: ``Kolmogorov
   must have arrived at his complexity definition before autumn
   1964, since Lyonya Bassalygo [Леонид Бассалыго] told me about
   it then. [Bassalygo confirms this; he remembers a walk during
   late autumn or early spring when Kolmogorov tried to explain
   him the complexity definition that was quite difficult to
   grasp at first.] On the other hand, it should be later than
   the randomness definition proposed in the Sankhya paper
   ~\cite{kolmogorov-1963} which was received April 1963 by the
   journal. Those considerations pin down the time of discovery
   to 1963--64, more exactly. (Kolmogorov never told me anything
   about the history of his discovery.)

   [On the other hand,] in his obituary note in the \emph{Journal
   of Applied Probability}, Vol.~25, No.~1, pp.~445--450, March
   1988, K.R.~Parthasarathy writes:

   ``Immediately after his arrival in Calcutta, Andrei
   Nikolaevich lost no time in plunging into discussions with
   the young students at the Institute about his recent research
   work on tables of random numbers, and the measurement of
   randomness of a sequence of numbers using ideas borrowed from
   mathematical logic. This piece of research was carried out by
   him during his travel by ship from the USSR to India; the
   ship was actually proceeding on an oceanographic
   expedition.''

   This seems to fix the time of the discovery of the complexity
   definition of randomness to 1962 [at least in some
   preliminary form] and to locate it to the ship that brought
   him to India for the reception of the degree of Doctor
   Honoris Causa at the University of Calcutta.''

   Kolmogorov gave several talks at the Moscow Mathematical
   Society but for most of them only the titles are known, and
   we may only guess what was there: \emph{Редукция данных с
   сохранением информации} (Data reduction that conserves
   information, March~22, 1961), \emph{Что такое
   ``информация''}? (What is information?, April~4, 1961),
   \emph{О таблицах случайных чисел} (On the tables of random
   numbers, October~24, 1962, probably corresponding to Sankhya
   paper~\cite{kolmogorov-1963}), \emph{Мера сложности конечных
   двоичных последовательностей} (A complexity measure for
   finite binary strings, April~24, 1963), \emph{Вычислимые
   функции и основания теории информации и теории вероятностей}
   (Computable functions and the foundations of information
   theory and probability theory, November~19, 1963),
   \emph{Асимптотика сложности конечных отрезков бесконечной
   последовательности} (Asymptotic behavior of the complexities
   of finite prefixes of an infinite sequence, December~15,
   1964; the title suggest that the last talk was about
   Martin-L\"of results, though Martin-L\"of remembers
   discussing these results with Kolmogorov
   only next spring, see below).
   Three later talks
   about algorithmic information theory (1968--1974)
   have short published
   abstracts (see Appendix~A.)}

The introduction of the complexity notion allowed to identify
randomness (for finite bit strings and fair coin) with
incompressibility. One should have in mind, however, that one
cannot hope to draw a sharp dividing line between random and
non-random strings of a given finite length, and the complexity
function $K(x)$ is defined up to a $O(1)$ term, so, strictly
speaking, only asymptotic statements are possible.

\section{Martin-L\"of definition of randomness}

To obtain such a sharp borderline one needs to consider infinite
sequences. A natural idea: to define randomness of an infinite
sequence in terms of complexity of its prefixes. The first
attempt was to say that a sequence $\omega_1\omega_2\ldots$ is
random if $K(\omega_1\ldots\omega_n)$ is maximal up to a
constant, i.e.,
        $$
K(\omega_1\ldots\omega_n)= n+O(1).
        $$
But Martin-L\"of\footnote{%
    Per Martin-L\"of, a mathematician from Sweden, studied
    Russian during his military service and then decided to make
    use of his knowledge by coming to Moscow and working with
    Kolmogorov.

    Martin-L\"of tells in~\cite{martin-lof-interview}: \ldots I
    had not worked on randomness before coming to Moscow in
    1964--65. Kolmogorov first gave me a statistical problem in
    discriminant analysis, which I solved, although I did not
    find it challenging enough. It was a problem that I
    might just as well have worked on at home in Stockholm. But
    I got to know Leonid (Lyonya) Bassalygo [Леонид Бассалыго],
    and he told me about Kolmogorov's new ideas about complexity
    and randomness, which I found very exciting. This was in
    late autumn 1964. So I started to learn the necessary
    recursive function theory from Uspenskij's
    book~\cite{uspensky-lectures}\ldots\ It was only when I told
    Kolmogorov about my first results on complexity oscillations
    in infinite binary sequences in early 1965 that complexity
    and randomness became the subject of our discussions. (So I
    did not learn about Kolmogorov complexity directly from
    Kolmogorov but only indirectly from Bassalygo).

    [As to the motivation,] I studied the previous literature on
    random sequences only after I had made my own first
    contributions. This resulted in the paper \emph{The
    Literature on von Mises' Kollektivs Revisited} published in
    the Swedish philosophical journal
    \emph{Theoria}~\cite{martin-lof-1969}. [As to the
    predecessors,] I have been most interested in Borel,
    particularly because he was the most important of the early
    French constructivists, which Brouwer called the
    pre-intuitionists. My affection for him may also have to do
    with the fact that I inherited a copy of Borel's
    \emph{Lecons sur la Th\'eorie des Fonctions}, with its many
    interesting Notes at the end, when my grandfather died in
    1958 and I was aged 16.

    When trying to require the complexities of the finite
    initial segments to be as big as possible, I discovered
    the unavoidable complexity oscillations about which I wrote
    my first paper on the subject (in Russian and typed by
    Nataliya Dmitrievna Svetlova [Наталья Дмитриевна Светлова
    (Солженицына)], who became Solzjenitsyn's wife in her second
    marriage). This led me to try the new approach of suitably
    interpreting the definition of null set in the sense of
    recursion theory. I should add that my primary reason for
    being interested in infinite rather than finite random
    sequences was to get rid of the additive constants that
    cropped up everywhere, and whose arbitrariness I found
    annoying. [This paper,] the first one of my two Russian
    papers was never published, but a typed copy of it should
    still exist somewhere in my unsorted archive. However, the
    results contained in it were subsequently published in
    English in my paper~\cite{martin-lof-1971}.

    The paper~\cite{martin-lof-1966r} is the second of the two
    papers that I have written in Russian. It summarizes a talk
    that I apparently gave in Moscow on 2~June~1965 and shows
    very clearly that I had not yet reached the definition of my
    Information and Control paper~\cite{martin-lof-definition} though
    I was on my way.

    Kolmogorov was immediately very interested in my two
    theorems on the unavoidable complexity oscillations in
    infinite binary sequences, which I told him about in the
    train on our way to Caucasus, more precisely, Bakuriani
    [Armenia] in early March~1965. In fact, he was so positive
    that he asked me to present my results as a sequel to a
    guest lecture that he gave in Tbilisi on our way back in
    late March. I do not think that he had thought himself about
    the problem of defining infinite random sequences by means
    of his complexity measure before then. So I think it is
    correct to say\ldots\ that he was more interested in finite
    random sequences. In a way, even if I have myself been
    interested in getting a good definition of randomness for
    infinite sequences, it is more striking that one can give a
    sensible definition of randomness already for finite
    sequences. Concerning finite random sequences, my own only
    contribution was the observation that the random elements of
    a finite population should be the ones whose conditional
    complexity given the population is maximal, that is,
    approximately equal to the logarithm to the base $2$ of the
    number of elements of the population, whereas Kolmogorov'
    original suggestion was to use the unconditional complexity.
    So, in the case of a completely random sequence of length
    $n$, we should use $K(x_1\ldots x_n|n)$ rather than
    $K(x_1\ldots x_n)$, and, in the case of Bernoulli sequences,
    $K(x_1\ldots x_n|n,s_n)$, where $s_n = x_1 +\ldots+ x_n$.

    I never had the opportunity of discussing my own definition of
    randomness for infinite sequences with Kolmogorov, simply
    because I did not find it until after I left Moscow in July
    1965. It must have been sometimes during the academic year
    1965--66. (End of quote.)
}
found that it is not possible (sequence with this property do
not exist).

Taking this difficulty into account, Martin-L\"of tried a
different approach and gave a definition of a random sequence
based on effectively null sets, making it more
measure-theoretic. The idea of this approach can be explained as
follows.

Let us define a random bit sequence (for simplicity we consider
only the case of a fair coin) as a sequence that satisfies all
probability laws. And probability law is a property of sequences
that is true for almost all sequences, i.e., for all sequences
outside some null set. Finally, a subset $X$ of the Cantor space
$\{0,1\}^\infty$ (of all infinite binary sequences) is a null set
if its uniform measure is $0$ (equivalent formulation: if for
every $\varepsilon>0$ there exists an infinite
sequence of intervals that
covers $X$ whose total measure is at most $\varepsilon$).

The problem with this definition is that random bit sequences
defined in this way do not exist at all. Indeed, for every
sequence $\alpha$ the singleton $\{\alpha\}$ is a null set, so
its complement $\{0,1\}^\infty\setminus\{\alpha\}$ can be
considered as a probability law, and $\alpha$ does not satisfy
this law.

Martin-L\"of pointed out that if we restrict ourselves to
\emph{effectively} null sets, this plan becomes quite
reasonable. A set $X$ is an effectively null set if there exists
an algorithm that (given positive rational $\varepsilon$)
generates a sequence of intervals that
cover $X$ and have total measure at most $\varepsilon$.
(Replacing algorithms with arbitrary functions, we get
a classical definition of null sets.)
It is easy to see that the union of all effectively null sets is
a null set, since there are only countably many algorithms.
Therefore random sequences (defined as sequences that do not
belong to any effectively null set) exist and the set of random
sequences has measure~$1$.

Moreover, Martin-L\"of have proved that the union of all
effectively null sets is an \emph{effectively} null set (in
other terms, there exists the largest effectively null set).
This maximal set consists of all non-random sequences. A set $X$
is effectively null if and only if $X$ is a subset of this
maximal effectively null set, i.e., $X$ does not contain any
random sequence.

We can formulate this in the following way. Let $P$ be some
property of binary sequences. Then the statements
     \begin{center}
$P(\alpha)$ is true for every random sequence $\alpha$
     \end{center}
and
     \begin{center}
the set of sequences $\alpha$ that do not satisfy $P$ is an effectively null set
     \end{center}
are equivalent in the word ``random'' in understood in Martin-L\"of sense.
This is nice because people often say, for
example, that ``for a random sequence $\alpha$ the limit
frequency is equal to $1/2$'' (the strong law of large numbers)
having in mind that the set of sequences that do not have this
property is a null set. Now this sentence can be understood
literally (if a null set is an effective null set, which is true
in most cases).

Martin-L\"of published this definition in 1966
in~\cite{martin-lof-definition}). His results were also
covered by a detailed survey paper~\cite{zvonkin-levin}.
written by two Kolmogorov's
young colleagues, Leonid Levin and Alexander Zvonkin (by Kolmogorov's
initiative; Kolmogorov carefully reviewed this paper once it
was finished and suggested many corrections).
This survey included
Martin-L\"of results as well as other results about complexity
and randomness obtained by the Kolmogorov school in Moscow. In
particular, a proof of the symmetry of information (an important
result obtained independently by Levin and Kolmogorov) was
included there.%
        \footnote{%
Levin recalls that being an undergraduate student he wanted to
convince Kolmogorov to be his
advisor and hoped that this result
would impress Kolmogorov. But Kolmogorov was rather busy, and
the appointment was postponed several times from February to
August 1967.
Finally, when Levin
called him again, Kolmogorov said
something like: ``O yes, come to see me, I have very interesting
results, the information is symmetric''. --- ``But, Andrei
Nikolaevich, this is exactly what I wanted to tell you.'' ---
``But do you know that the symmetry is only up to logarithmic
terms?'' --- ``Yes.'' --- ``And you can give a specific
example?'' --- ``Yes.'' Then Levin came to see Kolmogorov, they
discussed these results (later announced in
\cite{kolmogorov-1968-69} without proof; the first proof
appeared in~\cite{zvonkin-levin}). Levin indeed worked with
Kolmogorov during his undergraduate years and even earlier (the
first Levin's result was obtained under Kolmogorov's supervision
when Levin was in high school and published later
as~\cite{levin-1969}) but V.A.~Uspensky was officially listed
as his undergraduate advisor
for some formal reasons (see below).}

Martin-L\"of definition of randomness at first seems to be
purely measure-theoretic, it has nothing to do with selection
rules, martingales, and complexity. However, it turned out to be
closely related to these notions, and it was soon found by
different authors.

\section{Randomness and martingales: Schnorr}

During the next decade (1965--1975; recall that Kolmogorov
published his definition of complexity in 1965 and Martin-L\"of
published his definition of randomness in 1966) a lot of work
was done by different authors who provided missing links between
complexity, randomness and games (martingales). One of these
authors was C.P.~Schnorr.

As he tells~\cite{schnorr-talk}, after finishing his Ph.D. he
was looking for new topics. Martin-L\"of gave a course in
Erlangen, and the lecture notes of this course were distributed.
So this field become known in Germany, Schnorr heard a talk
about complexity and randomness and became interested. He wrote
several papers and then a book in Lecture Notes in Mathematics
series~\cite{schnorr-ln} based on his 1970 lectures (the book is
in German; it contains references to his other papers,
including~\cite{schnorr-unified} where many of the results from
the book are presented in English). His habilitation was based
on the results obtained in these papers.

In this book for the first time the notion of martingale was
used in connection with algorithmic randomness.\footnote{As
    Schnorr said in his talk~\cite{schnorr-talk}, he had not
    read Ville's book, but learned the notion of martingale
    indirectly through other sources.}
Schnorr defined a class of \emph{computable} (berechenbare) and
\emph{lower semicomputable} (subberechenbare) martingales. A
function $f$ (arguments are strings, values are reals) is called
computable if there is an algorithm that computes the values of
$f$ with any given precision: given $x$ and positive rational
$\varepsilon$, the algorithm computes some rational
$\varepsilon$-approximation to $f(x)$. A function is lower
semicomputable if there is an algorithm that, given $x$,
generates all rational numbers that are less than $x$. It is easy
to see that $f$ is computable if and only if both $f$ and $-f$
are lower semicomputable.

Schnorr then proved that a sequence is Martin-L\"of random if
and only if no semicomputable martingale wins against it, thus
providing a criterion of Martin-L\"of randomness in terms of
martingales. (A technical remark: note that the initial capital
can be non-computable in our setting.) Schnorr, however, was not
satisfied with this notion (lower semicomputability). He found
it rather counter-intuitive: there is no evident reason why we
should generate approximations from below (but not above) to
martingale values. So he thought that this class of martingales
is too broad and, therefore, the corresponding class of
sequences is too narrow. So he called Martin-L\"of random
sequences ``hyperzuf\"allig'' (``hyperrandom''; this name is not
in use now). He proved that there exists a sequence that wins
against all computable martingales but is not Martin-L\"of
random.

Schnorr also defined the notion of lower semicomputable
\emph{supermartingale}. A function $m$ is a supermartingale if
it satisfies the supermartingale inequality,
        $$
m(x)\ge \frac{m(x0)+m(x1)}{2}.
        $$
In game terms this means that player is allowed to throw away
her money during the game. Schnorr proved that lower
semicomputable supermartingales can be used for Martin-L\"of
randomness criterion in place of martingales.

Trying to find a better definition of randomness, Schnorr
considered a smaller class of effectively null sets (now called
sometimes ``Schnorr null sets'').
As we have
said, for an effectively null set $X$ there exists an algorithm
that given $\varepsilon>0$ generates a sequence of intervals
that cover $X$ and have total measure \emph{at
most}~$\varepsilon$. Schnorr introduced a stronger requirement:
this total measure should be \emph{equal to}~$\varepsilon$.
(This sounds a bit artificial; more natural equivalent
definition asks for a computably converging series of the length
of covering intervals.) The sequences that are outside all Schnorr null
sets are called ``zuf\"allig'' (now they are sometimes called
``Schnorr random'' sequences). Schnorr proved that this is
indeed a broader class of sequences than ``hyperzuf\"allig''
(Martin-L\"of random). He also proved a criterion in terms of
computable martingales: a sequence is zuf\"allig if and only if
no computable martingale ``computably wins'' on it (``computably
wins'' means that there exists a non-decreasing unbounded
computable function $h(n)$ such that the player's capital after
$n$ steps is greater than $h(n)$ for infinitely many $n$).

Schnorr's papers and book contain a lot of other interesting
things which were developed much later. For
example, he considers how fast player's capital increases during
the game and proves that if a sequence does not satisfy the
strong law of large numbers, then there exists a computable
martingale that wins exponentially fast against it (much later,
in 2000s, the growth of martingales was explored farther in
connection to the notions of effective dimension).

As Schnorr explains, one of his goals was to approach the notion
of ``pseudorandomness''. Sometimes even a sequence generated by
an algorithm looks similar to a random one; such sequences may
be used when the source of physical randomness is unavailable
and sometimes are called ``pseudorandom'', though this term may
have different more or less precise meanings. One of the
possible approaches to this phenomenon is that a
``pseudorandom'' object may have a short description, but the
time needed for the decompressing algorithm to process this
description is unreasonably large.\footnote{%
   Later a more practical theory of pseudorandom sequences was
   developed by Yao, Blum, Micali and others. Now it is a
   very important part
   of computational cryptography, see, e.g., the
   textbook~\cite{goldreich}. Schnorr later also worked in the
   field of computational cryptography.}
So Schnorr considers also complexity with bounded resources in
his book.

\section{Supermartingales and semimeasures}

Schnorr's lower semicomputable supermartingales are closely
related to other notion that appeared in Zvonkin and Levin's
1970 paper~\cite{zvonkin-levin}, the notion of a semicomputable
\emph{semimeasure}. It is easy to see that martingale (as
defined above) is just a ratio of two measures on the Cantor
space: an arbitrary one and the uniform one. More formally, let
$Q$ be any measure on Cantor space and let $P$ be the uniform
Bernoulli measure. Then the ratio $Q(I_x)/P(I_x)$, where $I_x$
is the interval rooted at binary string $x$ (the set of all
extensions of $x$), is a martingale. Moreover, every martingale
can be represented in this way. The supermartingales correspond
in the same way to objects that Levin called ``semimeasures''.

A semimeasure is a measure on the set $\Sigma$
of all finite and infinite binary sequences. Let $\Sigma_x$ be
the set of all extensions (finite and infinite) of a binary
string $x$. Then $\Sigma_x=\Sigma_{x0}\cup
\Sigma_{x1}\cup\{x\}$. If $Q$ is a measure on $\Sigma$, the
inequality
        $$
Q(\Sigma_x)\ge Q(\Sigma_{x0})+Q(\Sigma_{x1})
        $$
holds; moreover, any non-negative real-valued function $q$ on
finite strings that satisfies the inequality $q(x)\ge
q(x0)+q(x1)$, determines a measure on~$\Sigma$. The
difference between both sides of this inequality is the
probability of the finite string $x$. A semimeasure is
\emph{lower semicomputable} if the function $x\mapsto
q(x)=Q(\Sigma_x)$ is lower semicomputable.

Lower semicomputable semimeasures are considered
in~\cite{zvonkin-levin}; Levin proved that they can be
equivalently defined as output distributions of probabilistic
machines that have no input, use internal fair coin and generate
their output sequentially (bit by bit). Levin proved also that
there exists a maximal lower semicomputable semimeasure
(\emph{universal semimeasure}, sometimes called \emph{a priori
probability} on the binary tree). This notion can be also
considered as a formalization of Solomonoff's ideas.

The connection between semimeasures and supermartingales:
supermartingales can be defined as fractions where the numerator
is a semimeasure and denominator is the uniform Bernoulli
measure (similar to the description of martingales as fractions
of two measures). Lower semicomputable semimeasures correspond
to lower semicomputable supermartingales. This representation of
(semi)martingales as ratios can be easily generalized to other
probability distributions, e.g., to the case of a non-symmetric
coin. If $P$ is the distribution declared by the game organizers
(now not necessarily uniform), then in the ``fair'' game the
player's capital is a $P$-martingale, i.e., the ratio $Q/P$
where $Q$ is some measure. (The notion of martingale with
respect to a non-uniform measure was also considered by Schnorr
in~\cite{schnorr-ln}.)

In a similar way $P$-\emph{super}martingales (that allow the
player to discard some money at each step) can be defined as
ratios $Q/P$ where $Q$ is a \emph{semi}measure. This implies,
for example, that for any computable measure $P$ there exists a
maximal lower semicomputable $P$-supermartingale: it is the
ratio $A/P$ where $A$ is the a priori probability (the largest
lower semicomputable semimeasure). The last observation provides
a connection between maximal $P$-supermartingales for different
$P$; as Levin points in one of the letters to Kolmogorov (see
the Appendix) the advantage of the a priori probability notion
is that the same notion can be compared to different measures.
When switching from (semi)measures to (super)martingales one
object (the a priori probability) is transformed into a family
of seemingly different objects (maximal lower semicomputable
supermartingales with respect to different computable measures).

However, a natural goal: ``to obtain
a criterion of randomness (for
infinite sequences) in terms of complexity of their prefixes''
(the idea to relate complexity and randomness
was present already in the 1965 Kolmogorov
publication~\cite{kolmogorov-1965})
was not achieved either in Zvonkin and Levin paper
or in Schnorr's book. This was done few years later when new
versions of complexity (monotone and prefix complexities)
appeared.

\section{Prefix complexity}

Prefix complexity was introduced by Levin and Chaitin. Since the
introduction of prefix complexity sometimes becomes a source of
unnecessary controversy, some historical
clarifications would be useful here. To put the story short, the
first publications where (1)~the prefix complexity was defined
in terms of self-delimiting codes and as the logarithm of the
maximal lower semicomputable converging series, and (2)~the
claim that these definitions coincide was made (without proofs),
are~\cite{levin-1974,gacs-1974}. These publications appeared in
1974 in
Russian; English translations of these two papers were published
in 1976 and 1975 respectively (see~\cite{gacs-review}); the
logarithm of the maximal lower semicomputable
converging series
(but not the
self-delimiting descriptions) was considered also in unpublished
thesis of Levin in 1971.\footnote{%
   Let us add some historical remarks about situation in the
   Mathematics Department of Moscow State University and in
   Russia at that time. The typical track of a future
   mathematician at that time was 5 years of undergraduate
   studies (высшее образование) plus 3 years of graduate school
   (аспирантура). After the graduate school student is assumed
   to defend a thesis and get a title ``kandidat
   fiziko-matematicheskih nauk'' (кандидат физико-математических
   наук) which is a rough equivalent of Ph.D. Unlike the US
   universities, the student of Moscow State University (and
   other Soviet universities) had to decide what is his major
   before entering the university: e.g., the mathematics and
   physics programs are administered by different departments,
   have no common courses, different entrance procedures etc.
   After two years of undergraduate studies at mathematics
   department, a student had to choose a division (кафедра)
   which he wants to join for three remaining years, and a
   scientific advisor in the chosen division. (It could be, say,
   Algebra Division, or Geometry and Topology Division, etc.) At
   the end of the 5th year student writes a thesis (дипломная
   работа). Sometimes this thesis is considered as something
   close to the Master thesis in the US.

   To enter the graduate school after finishing 5 years of
   undergraduate studies, one needed a good academic record and
   (a very important condition!) a recommendation from the local
   communist party and komsomol (комсомол) organization.
   Komsomol (an abbreviation for \textbf{ком}мунистический
   \textbf{со}юз \textbf{мол}одёжи, communist union of the young
   people), like Hitlerjugend in Germany, was
   almost obligatory, and included people of age 14--28, so most
   university students were komsomol members (комсомольцы),
   though there were some exceptions and this requirement was
   never formalized as a law.

   Levin was a student of a special boarding school founded by
   Kolmogorov (unofficially called Kolmogorov's boarding school,
   колмогоровский интернат); during 1963/4 academic year he was
   a student of a similar school in Kiev (now Ukraine) and then
   managed to move to Moscow for 1964/5 academic year. Then (in
   January 1966) he entered the Moscow university becoming a
   first-year undergraduate in the middle of the academic year
   (there was some exceptional procedure for the students of
   Kolmogorov's school in this year related to the change in the
   education system in the USSR that moved from 11-years to
   10-years education program).

   Being not only Jewish (already a handicap at that
   time) but also a kind of non-conformist, Levin as an
   undergraduate student created a lot of troubles for local
   university authorities. As a member of komsomol, he became
   elected local komsomol leader but he defied the policies
   established by
   the Communist Party supervisors (and this was mentioned
   in his graduation letter of recommendation, a very important
   document at the time).
   No wonder he was effectively   %! ``effectively'' added
   barred from applying to any
   graduate school when he finished undergraduate studies
   at the Mathematical Logic Division (кафедра математической
   логики) in 1970. (His official undergraduate advisor was Vladimir
   A.~Uspensky who was Kolmogorov's student in 1950s.
   Kolmogorov officially
   did not belong to Mathematical Logic division and asked his
   former student Uspensky to replace him in this capacity.)
   However, Kolmogorov managed to secure a research scientist
   position for Levin
   (with the help of the University rector, a prominent
   mathematician and a very decent person, I.G.~Petrovsky) in
   the University statistical laboratory (Kolmogorov was a head
   of this laboratory).

   Being there, in 1971 Levin wrote a ``kandidat''
   thesis (that contained
   mostly Levin's results included in~\cite{zvonkin-levin}, but
   also some others, including the probabilistic definition of
   prefix complexity) and tried to find a place for its defense.
   (According to the rules, the thesis defence was not technically
   connected to a graduate school (if any) of defendant's affiliation,
   only a recommendation from the institution where dissertation
   was prepared was required;
   in this case the person was called ``соискатель''. Though most
   graduate students in the USSR were
   defending their thesis in the same institution
   (sometimes a few years later after their studies in the graduate
   school), the thesis defense was not a university
   affair, but regulated by a special government institution,
   called ``Высшая Аттестационная Комиссия''.)   %! quotes added, ВАК deleted

   In Moscow it was clearly impossible, and finally the defense
   took place in Novosibirsk (in Siberia).
   The thesis received strong approvals from official reviewers
   (J.~Barzdin, B.~Trachtenbrot and his lab), the reviewing
   institutions (Leningrad Division of Steklov Mathematical
   Institute) and the advisor (Kolmogorov and his lab).
   Nevertheless, the defence
   was unsuccessful
   (quite untypical event). According to Levin, the most active
   negative role during the council meeting was played by
   Yu.L.~Ershov (recursion theorist, now a member of the Russian
   Academy of Sciences) but Levin believes that Ershov did not
   have other choice unless he was ready to get into career
   troubles himself; however, Ershov did also something ``above
   and beyond the call of duty'' (as Levin puts it) as a Soviet
   scientific functionaire~--- he insisted that the ``unclear
   political position'' of Levin should be mentioned in the
   council decision. This effectively prevented Levin's
   defense in any other place in the Soviet Union (even with a
   new thesis) and therefore barred a scientific career in
   Soviet Union for him. Fortunately, Levin got a permission to
   leave Soviet Union and emigrated to US where he got many well
   known results in different areas of theoretical computer
   science (about one-way functions, holographic proofs et al.).
   As Levin recalls, KGB made it
   known that they think going away would be the best option for
   him; they even asked Kolmogorov to deliver this advice (which
   Kolmogorov did, though he did not indicated whether he himself
   agrees\ldots)
   Now we can make jokes about these events (Levin once noted
   that a posteriori Ershov's behaviour was a favor for him: it
   was a motivation to leave Soviet Union) but at that time
   things were much more dramatic.

   But while being still in the USSR after this unsuccessful
   defense, Levin followed an advice of some friend, who told
   that Levin should publish his results while he is still
   allowed to publish papers in Soviet journals (this was not a
   joke, the danger was quite real) and published a bunch of
   papers in 1973--1977. These papers were rather short and
   cryptic, a lot of things was stated there without proofs, so
   many ideas from them were really understood only much later.

   We will trace only two main contributions made in these
   papers: the prefix complexity, and the randomness criterion
   in terms of monotone complexity.}
In 1970~paper~\cite{zvonkin-levin}
an a priori probability (on a binary tree, as defined in this
paper) of a sequence $0^n1$ is considered (last paragraph on
p.~107) and some properties of this quantity
are proved, though no name is given
for it; this quantity coincides with a maximal lower semicomputable
converging series (up to $O(1)$ factor, as usual).

At the same time Chaitin independently came to the same two
definitions (self-delimited complexity and logarithm of
probability) in~\cite{chaitin-1975}; this paper, submitted in 1974,
contained, among other results, the first published proof of
their equivalence. (See more about the history of this paper
below.)

The prefix complexity, as we have said,
can be defined in different ways. The first approach defines
prefix complexity of $x$ as the length of the shortest program
that produces $x$, but the programming language must satisfy an
additional requirement. In Levin's paper~\cite{levin-1974} this
requirement is formulated as follows: if a bit string $p$
considered as a program produces some output $x$, then its
extensions either produce the same $x$ or do not produce
anything. The 1974 paper refers for details to Gacs' paper
of the same year~\cite{gacs-1974}\footnote{%
        Peter Gacs came to Moscow State University for 1972/3
        academic year from Hungary where he became interested in
        this topic after reading
        Kolmogorov paper~\cite{kolmogorov-1965}, Martin-L\"of
        lecture notes from Erlangen and Zvonkin and Levin's
        paper~\cite{zvonkin-levin} and started correspondence
        with Levin by sending him some paper about randomness
        characterization in terms of complexity. When Gacs came
        to Moscow in 1972, Levin explained his criterion of
        randomness in terms of monotone complexity which looked
        much better to Gacs so his paper was never published.
        Then Levin explained the notion of prefix complexity to
        Gacs and asked whether it is symmetric (with $O(1)$
        precision). The negative answer obtained by Gacs became
        part of his paper~\cite{gacs-1974} that included also
        some Levin's results, including the $O(1)$-formula for the
        prefix complexity of a pair (attributed to
        Levin). The prefix complexity is very briefly introduced
        in the beginning of this paper with the remark
        ``considered in detail by Levin''.}
and to other Levin's paper (then
unpublished; it was published only in 1976~\cite{levin-1976}).
In Chaitin's paper mentioned above~\cite{chaitin-1975}\footnote{%
        This paper was written~\cite{chaitin-home} in 1974
        during the visit to the IBM Watson Lab in Yorktown
        Heights for a few months. Chaitin's work there has
        another important implication: an unpublished manuscript
        by R.~Solovay~\cite{solovay-manuscript}. In his
        talk~\cite{calude-talk} Cristian Calude tells the story:
        ``When I started reading and trying to understand the
        subject to write my book ``Information and
        Randomness''~\cite{calude-book}, I discussed this with
        Greg Chaitin and he told me: look, if you really want to
        write a good book, you have to read Solovay's
        manuscript\ldots\ So I started asking around, and
        eventually wrote to Solovay: Greg Chaitin told me that I
        should read your manuscript; could I have a copy?
        Solovay answered: I had one, but I don't have it
        any more. This was in 1991, I think. I tried again to
        get it and eventually I contacted Charles Bennett, and
        he had one copy; he was very kind to send me a copy of
        this copy. That is also an interesting story which Greg
        Chaitin told me about how this book [manuscript] was
        written. Solovay was for one year at IBM on a sabbatical
        leave and he was asked to write a report about Chaitin's
        work. Probably most of us would write a report of
        two or three pages
        and forget forever about it. But Solovay
        took it very seriously, so he rewrote many parts of the
        theory in his completely different new style, and he
        solved also a substantial number of open problems at
        that stage. This was a kind of shock: look, this guy
        is so bright, he has nothing to do with this field, he
        comes, he reads this bunch of papers, he produces this
        beautiful manuscript solving so many problems and at the
        end of the day he does not want even to publish
        anything! Solovay never
        published this manuscript. I sent Solovay a copy of his `lost'
        manuscript
        and he said: well, if you have a student or whoever
        would like to read and edit and publish the book, fine
        with me, but I am not interested in working on it. It
        had to wait until Rod Downey and Denis Hirschfeldt had
        the force to get through and recuperate most of the
        results in this manuscript.''}
a slightly different requirement is used: if a bit string $p$
considered as a program outputs $x$, then none of $p$'s
extension could produce any output. Both restrictions reflect
the intuitive idea of a self-delimiting program (that does not
contain an end-marker; the machine should be able to find out
when the program ends) though in technically different ways.

Another way to define prefix complexity uses probabilities; as
we have mentioned, it appeared in Levin's thesis~(1971) that
remained unpublished. Consider the lower semicomputable series
of non-negative reals with sum at most~$1$ ($\sum p_n\le 1$
where $p_n\ge 0$ and the function $n\mapsto p_n$ is lower
semicomputable). These series correspond to machines that use
internal fair coin to produce some integer (or, may be, do      %! $n$ deleted
not produce anything) if we let $p_n$ be the probability of
output $n$.

Among those series there exists a maximal one (up to $O(1)$
factor). It is called a priori probability on integers (and is
closely related to the a priori probability on bit strings
considered in Zvonkin--Levin paper~\cite{zvonkin-levin}: a
priori probability of a bit string $0^n1$ coincides with the a
priori probability of integer~$n$ up to $O(1)$ factor).

A very important property of these notions: minus binary
logarithm of an a priori probability equals prefix complexity
(up to $O(1)$ additive term). This property is mentioned without
proof both in~\cite{gacs-1974} and~\cite{levin-1974}; the proof
was published for the first time in~\cite{chaitin-1975}. This
proof implies also that two version of prefix-free requirements
mentioned above lead to the same complexity function (up to
$O(1)$ additive term).

Another advantage of prefix complexity, also discovered
independently by Levin (the proof, attributed to Levin, is
published in~\cite{gacs-1974}) and Chaitin (the proof is
published in~\cite{chaitin-1975}) is a more precise (up to
$O(1)$-term) formula for the complexity of a pair in terms of
conditional complexities. This formula is an improvement of the
symmetry of information theorem that was earlier proved for
plain complexity with bigger (logarithmic) error terms by
Kolmogorov and Levin.

\section{Randomness criterion: Schnorr and Levin}

It was soon understood by Schnorr and Levin that the original
goal of describing randomness in terms of complexity can be
achieved if one changes a bit the definition of complexity
making it monotonic in some sense.

Schnorr suggested such a modification in his talk at 4th STOC in
1972~\cite{schnorr-1972}. The idea of the modification was to
take into account that prefixes of a sequence are not separate
binary strings but prefixes of one infinite sequence. As Schnorr
puts it (\cite{schnorr-1972}, pp.~168--169), ``it has already
been observed that there must be some difference in the concept
of regularity of finite objects which do not involve a direction
(for instance a natural number) and the concept of regularity of
infinite sequences (as well as finite subsequences [prefixes] of
an infinite sequence) where a natural direction is involved. For
example, he who wants to understand a book will not read it
backwards, since the comments or facts which are given in his
first part will help him to understand subsequent chapters (this
means they help him to find regularities in the rest of the
book). Hence anyone who tries to detect regularities in a
process (for example an infinite sequences or an extremely long
finite sequence) proceeds in the direction of the process.
Regularities that have ever been found in an initial segment of
the process are regularities for ever. Our main argument is that
the interpretation of a process (for example to measure his
complexity) is a process itself that proceeds in the same
direction.''\footnote{%
        This argument sounds convincing; however, one may expect
        that randomness of a binary sequence is invariant under
        computable permutation of its terms while Schnorr's
        criterion of randomness in terms of monotone complexity
        is not. Recently A.~Rumyantsev pointed out the following
        simple invariant criterion: $\textit{KP}(A,\omega(A))\ge
        |A|-O(1)$. Here $\textit{KP}$ stands for the prefix complexity
        of a pair; $A$ is a finite set of indices of size $|A|$
        and $\omega(A)$ is a restriction of $\omega$ onto $A$ (a
        bit string of length $|A|$).}
Then he gives a formal definition of monotone
complexity, called ``process complexity'' in his paper, and
notes that ``basic properties of processes have been developed
independently in [5] and [8]'' (i.e., \cite{schnorr-ln}
and~\cite{zvonkin-levin} in our list; note that none of these two
publications includes a definition of monotone/process complexity).

Using his definition, Schnorr proves that a sequence in
Martin-L\"of random if and only if its $n$-bit prefix has
monotone complexity $n+O(1)$.

Levin~\cite{levin-monotone-1973} proves essentially the same
result using a slightly different version of the monotone
complexity (used also in subsequent paper of
Schnorr~\cite{schnorr-1977}). Levin also notes that the same
proof works for the so-called ``a priori complexity'', the minus
logarithm of the a priori probability on the binary tree. This
statement is equivalent to Schnorr's characterization of
randomness in terms of semicomputable supermartingales (though
Levin does not say anything about martingales).

Chaitin in~\cite{chaitin-1975} suggested prefix complexity as a
tool to define randomness. He calls an infinite sequence
$\omega_1\omega_2\ldots$ random if there exists $c$ such that
        $$
H(\omega_1\ldots\omega_n)\ge n-c
        $$
for all $n$ (he used letter $H$ to denote prefix complexity;
Levin used $KP$; now the letter $K$ is most often used), and
writes: ``C.P.~Schnorr (private communication) has shown that
this complexity-based definition of a random infinite string and
P.~Martin-L\"of statistical definition of this concept are
equivalent''. As Schnorr remembers in his
talk~\cite{schnorr-talk}, ``I knew the first paper of Chaitin
that has been published one year later after the Kolmogorov's
1965 paper but it was the next paper which really made Chaitin
also one of the basic investigators of complexity. This was a
paper on self-delimiting or prefix-free descriptions and this
was published in 1975 in the Journal of the ACM. In fact I was a
referee of this paper and I think Chaitin knew this because I've
sent my personal comments and suggestions to him and he used
them''.

\section{Lower semicomputable random reals}

One more result about randomness in~\cite{chaitin-1975} is an
example of a lower semicomputable random real number, now well
known as ``Chaitin's $\Omega$ number''. It is related to a
philosophical question: can we specify somehow an individual
random sequence? One would expect at first the negative answer:
if a sequence has some description that defines it uniquely, how
can we treat it as random?

This negative answer is supported by the (evident) result: a
computable sequence is not Martin-L\"of random (for the case of
a fair coin, i.e., the uniform Bernoulli distribution). However,
if we do not insist that description is an algorithm that
computes our sequence and let it be less direct, the answer
becomes positive. Indeed, in~\cite{zvonkin-levin} the following
result attributed to Martin-L\"of is stated (Theorem~4.5):
there exists a $\Sigma^0_2$-sequence that is Martin-L\"of
random. This means that there exist a decidable property
$R(n,p,q)$ of three natural numbers such that the sequence
$\omega$ defined by equivalence
        $$
\omega_n = 1 \ \Leftrightarrow\ \exists p\,\forall q\, R(n,p,q)
        $$
is Martin-L\"of random. This provides an example of an
individual explicitly described (though in a non-constructive
way) random sequence.

The example of a random $\Sigma_2^0$-sequence appears also in
Theorem~4.3 in Chaitin's 1975~paper~\cite{chaitin-1975}, but
Chaitin went farther in this direction. He noticed that a
Martin-L\"of random sequence can be a binary representation of a
lower semicomputable real number. Speaking about random reals,
we identify real numbers in the interval $(0,1)$ with their binary
representations. (The
collisions like $0.0011111\ldots=0.0100000\ldots$ do not matter
since this can happen only for non-random sequences.) Recall
that a real number $x$ is lower semicomputable if there is an
algorithm that enumerates all rational numbers less than~$x$.
(Equivalent definition: if $x$ is a limit of an increasing
computable sequence of rational numbers.) It is easy to see that
all lower semicomputable reals $x\in(0,1)$ have binary
representations in $\Sigma^0_2$ but the reverse statement is not
true.

        %! next paragraph changed to describe the situation more exactly
This alone wouldn't make Chaitin's example of lower
semicomputable random real so popular. In fact, Section 4.4 of
\cite{zvonkin-levin} (proof ot Theorem~4.5 mentioned above)
already
constructs a specific example of a
random real, i.e., the smallest real outside an
effective open set of small measure that covers all non-random
reals. Zvonkin and Levin used the language of binary
sequences, not reals (which makes the description a bit more tedious)
and did not mention explicitly the
lower semicomputability (which follows immediately from the construction).
But the main reason why Chaitin's example became
so famous is in the form of the description. Chaitin's lower
semicomputable real $\Omega$ has simple and intuitive meaning:
it is the probability that the universal machine used in the
definition of prefix complexity terminates on a randomly chosen
program. This could create an impression
that we really have a
random real ``in our hands'': this is the probability of the
event ``the universal machine terminates on random
input''.\footnote{%
        A similar thing was done once to test early Unix
        utilities: they were fed with random bits and
        crashed quite often! In fact, standard programming
        languages and executable file formats satisfy Chaitin's
        requirements for universal machine if we
        ignore that machine word has finite size, usually between $8$
        and~$64$~bits.}

\section{Subsequent achievements}

The study of randomness as a mathematical object had clearly a
philosophical motivation related to the foundations of
probability theory. However, the mathematical theory has its own
logic of development: answering some philosophically motivated
questions, it introduces new notions and new questions related
to these notions. So the mathematical theory of randomness (and
related algorithmic information theory) became a rich
mathematical subject. In the last decade it attracted a lot of
attention from the recursion theorists who used advanced
techniques developed in recursion theory to understand the
randomness definitions better. For example, they looked at one
of the first definitions of randomness (from Kolmogorov's
papers) and proved that it coincides with Martin-L\"of
randomness relativized to
$\mathbf{0}'$-oracle~\cite{nies,miller}.

The other thread that has some philosophical and historical
interest is related to non-monotonic selection rules and
martingales. In Mises definition the terms of the sequence are
revealed in some fixed order (time order, if we look at casino's
example). He never explicitly mentioned other possibilities
(though he sometimes writes about data whose ordering is not
clear, like statistical data about deaths used by an
insurance company). When he was forced to provide a formal
definition of a selection rule, this monotonicity is explicitly
present in the definition.

However, one can consider other examples that motivate
non-monotonic selection. Imagine that casino prepares random
bits and write them on cards which are then placed on a table
(so that bits are invisible). The player is then allowed to look
at the cards in any order and also make bets (before the card is
turned). Imagine that she manages to win systematically; does it
implies that the sequence is not random?

As D.~Loveland~\cite{loveland-1966} explains this: ``Consider
the following ``practical'' situation. A manufacturer produces
very cheaply and quickly some item which has a large fluctuation
in life expectancy from item to item, with the fluctuation
passing through a threshold of acceptance. The producer would
naturally wish to cull out the unaccepted items but (it is
presumed) cannot test the item to be used for life expectancy
without destroying it. He must then look for ``systematic
fluctuations'' in the process so that he can select the items to
be used based on the knowledge of the process including
knowledge of tested items then ineligible for use. If the
process were random in the aforementioned sense, then no system
of testing previously manufactured items would indicate whether
the next item manufactured should be chosen for use or whether
one should choose, rather, some future item after more testing.
However, suppose the manufacturer numbers each item
consecutively as it is produced and allows it to fall it into a
bin from which items are drawn to be tested or selected for use.
Then he may test higher numbered items before digging down in
the bin to select a specific item for use.''

Earlier the same extension was suggested by Kolmogorov in a
footnote in his paper~\cite{kolmogorov-1963}. It led to many
interesting questions. For example, how complex should be
prefixes of a sequence that is random in the sense of
Mises--Church definition and in this extended
Mises--Kolmogorov--Loveland definition? Kolmogorov
claimed~\cite{kolmogorov-1968-69} that in both cases complexity
could be logarithmic, but later An.~Muchnik has shown that it is
not the case (see~\cite{uspensky-semenov-shen}) for Mises--Kolmogorov--Loveland
randomness (while for Mises--Church randomness Kolmogorov was right).
        %! adding explanation what Muchnik did

% Laurent: future developments?

Many other interesting results are obtained but their
description goes far beyond the scope of this paper.

\section{Concluding remarks}

Remember that Mises' initial reason to consider collectives was
the desire to explain what probability is and why and how the
mathematical probability theory can be applied to the real
world. The question ``why'' is rather philosophical one, but one
can try to answer to second part, ``how'', and describe the
current best practice. Here is an attempt to provide such a
description taken from~\cite{uspensky-semenov-shen,shen}.

``The application of probability theory has two stages. At the
first stage we try to estimate the concordance between some
statistical hypothesis and experimental results. The rule ``the
actual occurrence of an event to which a certain statistical
hypothesis attributes a small probability is an argument against
this hypothesis'' (Polya~\cite{polya}, Vol.~II, Ch.~XIV, part~7,
p.~76), it seems, could be made more correct if we are allowed
to consider only ``simply described'' events. It is clear that
the event ``$1000$ tails appeared'' can be described more simply
that the event ``a sequence $A$ appeared'' where $A$ is a
``random'' sequence of $1000$ heads and tails (these two events
have the same probability). This difference may explain why our
reactions to these events (we have in mind the hypothesis of a
fair coin) are so different. To clarify the notion of a ``simply
described event'' the notion of complexity of the constructive
object (introduced by Kolmogorov) may be useful.

Let us assume that we have already chosen a statistical
hypothesis concordant (as we think) with the result of
observations. Then we come to the second stage and derive some
conclusions from the hypothesis chosen. Here we have to admit
that probability theory makes no predictions but can only
recommend something: if the probability (computed on the basis
of the statistical hypothesis) or an event $A$ is greater than
the probability of an event~$B$, then the possibility of the
event~$A$ must be taken into consideration to a greater extent
than the possibility of the event~$B$.

One can conclude that events with very small probabilities may
be ignored. Borel~\cite{borel-1913} writes ``\ldots Fewer than a
million people live in Paris. Newspapers daily inform us about
the strange events or accidents that happen to some of them. Our
life would be impossible if we were afraid of all adventures we
read about. So one can say that from a practical viewpoint we
can ignore events with probability less that one
millionth\ldots\ Often trying to avoid something bad we are
confronted with even worse\ldots\ To avoid this we must know
well the probabilities of different events'' (Russian ed.,
pp.~159--160).

% Laurent: it is somehow shameful to cite Borel using English
% translation of Russian edition -- may be you have the original
% reference?

Sometimes the criterion for selection of a statistical
hypothesis and the rule for its application are united in the
statement ``events with small probabilities do not happen''. For
example, Borel writes ``One must not be afraid to use the word
``certainty'' to designate a probability that is sufficiently
close to $1$.'' (\cite{borel-1950}, Russian ed., p.~7). But we
prefer to distinguish between these two stages, because at the
first stage the existence of a simple description of an event
with small probability is important, and at the second stated it
seems unimportant. (We can expect, however, that events
interesting to us have simple descriptions because of their
interest.)''

This description (which, we believe, still describes adequately
the current best practice of probability theory application)
uses the notions of algorithmic information theory only once
(when describing when we reject a statistical hypothesis), but
this use seems to be important.

Let us note also that this description shows that quantum
mechanics does not make a real difference compared to
probability theory and statistical mechanics: we just replace
``small probability'' by ``small amplitude'' in the scheme
described. (However, to provide a foundation for the measurement
procedure, one should prove a quantum counterpart for the law of
large numbers: the amplitude of the event ``measured frequency
of some outcome diverges significantly from the square of the
assumed amplitude of this outcome'' is small.)

More detailed discussion can be found in
\cite{shen-reachability}.

\section*{Appendix A: Abstracts of Kolmogorov's talks}

Some talks at the meetings of Moscow Mathematical Society have
short abstracts published in the journal ``Успехи математических
наук'' (Uspekhi matematicheckikh nauk, partially translated as
``Russian mathemathical surveys''; these abstracts were not
translated). Here we reproduce abstracts of three talks given by
A.N.~Kolmogorov devoted to algorithmic information theory
(translated by Leonid Levin).

\bigskip

\textbf{I.}
[vol.~23, no.~2, March-April 1968].

\medskip

1. A.N.~Kolmogorov, ``Several theorems about algorithmic entropy
and algorithmic amount of information''.

Algorithmic approach to the foundations of information theory
and probability theory was not developed far in several years
from its appearance since some questions raised at the very
start remained unanswered. Now the situation has changed
somewhat. In particular, it is ascertained that the
decomposition of entropy $H(x,y)\sim H(x)+H(y|x)$ and the
formula $J(x|y)\sim J(y|x)$ hold in algorithmic concept only
with accuracy $O([\log H(x,y)])$ (Levin, Kolmogorov).

Stated earlier cardinal distinction of algorithmic definition of
a Bernoulli sequence (a simplest collective) from the definition
of Mises-Church is concretized in the form of a theorem: there
exist Bernoulli (in the sense of Mises-Church) sequences
$x=(x_1,x_2,...)$ with density of ones $p=\frac 1 2$, with
initial segments of entropy (``complexity'')
$H(x^n)=H(x_1,x_2,...,x_n)=O(\log n)$ (Kolmogorov).

For understanding of the talk an intuitive, not formal,
familiarity with the concept of a computable function suffices.

(Moscow Mathematical Society meeting, October~31, 1967)

\bigskip

\textbf{II.}
[vol.~27, no.~2, 1972]

\medskip

{\bf 1}. A.N. Kolmogorov. ``Complexity of specifying and complexity of
constructing mathematical objects''.

\begin{enumerate}

\item Organizing machine computations requires dealing with
 evaluation of (a)~complexity of programs, (b)~the size of
 memory used, (c)~duration of computation. The talk describes a
 group of works that consider similar concepts in a more
 abstract manner.

\item It was noticed in 1964-1965 that the minimal length $K(x)$
 of binary representation of a program specifying construction
 of an object $x$ can be defined invariantly up to an additive
 constant (Solomonoff, A.N.~Kolmogorov). This permitted using
 the concept of {\em definition complexity} $K(x)$ of
 constructive mathematical objects as a base for a new approach
 to foundations of information theory (A.N. Kolmogorov, Levin)
 and probability theory (A.N. Kolmogorov, Martin-L\"of, Schnorr,
 Levin).

\item Such characteristics as ``required memory volume,'' or
 ``required duration of work'' are harder to free of technical
 peculiarities of special machine types. But some results may
 already be extracted from axiomatic ``machine-independent''
 theory of broad class of similar characteristics (Blum, 1967).
 Let $\Pi(p)$ be a characteristic of ``construction complexity''
 of the object $x=A(p)$ by a program $p$, and $\Lambda(p)$
 denotes the length of program $p$. The formula $K^n\Pi(x)=
 \inf(\Lambda(p): x=A(p), \Pi(p)=n)$ defines ``$n$-complexity of
 definition'' of object $x$ (for unsatisfiable condition the
 $\inf$ is considered infinite).

\item Barzdin's Theorem on the complexity $K(M_\alpha)$ of prefixes
 $M_\alpha$ of an enumerable set of natural numbers (1968) and
 results of Barzdin, Kanovich, and Petri on corresponding
 complexities $K^n\Pi(M_\alpha)$, are of general mathematical
 interest, as they shed some new light on the role of extending
 previously used formalizations in the development of mathematics.
 The survey of the state of this circle of problems was given in
 the form free from cumbersome technical apparatus.

\end{enumerate}

(Moscow Mathematical Society meeting, November 23, 1971)

\bigskip
\textbf{III.}
[Vol.~29,. no.~4 (155), 1974]
\medskip

 1. A.N. Kolmogorov.
 ``Complexity of algorithms and objective definition of randomness''.

To each constructive object corresponds a function $\Phi_x(k)$
of a natural number $k$ -- the log of minimal cardinality of
$x$-containing sets that allow definitions of complexity at most
$k$. If the element $x$ itself allows a simple definition, then
the function $\Phi$ drops to $1$ even for small $k$. Lacking
such a definition, the element is ``random'' in a negative
sense. But it is positively ``probabilistically random'' only
when function $\Phi$, having taken the value $\Phi_0$ at a
relatively small $k=k_0$, then changes approximately as
$\Phi(k)=\Phi_0-(k-k_0)$.

(Moscow Mathematical Society meeting, April~16, 1974)

\section*{Appendix B. Levin's letters to Kolmogorov}

These letters do not have dates but were written after
submission of~\cite{zvonkin-levin} in August 1970 and
before
Kolmogorov went (in January~1971)
to the oceanographic expedition (``Dmitry
Mendeleev'' ship). Copies provided by L.~Levin
(and translated by A.~Shen).

\subsection*{I.}

Dear Andrei Nikolaevich! Few days ago I've obtained a result
that I like a lot. May be it could be useful to you if you
work on these topics while traveling on the ship.

This result gives a
formulation for the foundations of probability theory different from Martin-L\"of. I think
it is closer to your initial idea about the relation between
complexity and randomness and is much clearer from the
philosophical point of view (as, e.g., [Yu.~T.]~Medvedev says).

Martin-L\"of considered (for an arbitrary computable
measure~$P$) an algorithm that studies a given sequence and finds
more and more deviation from $P$-randomness hypothesis. Such an
algorithm should be $P$-consistent, i.e., find deviations of
size $m$ only for sequences in a set that has measure at most
$2^{-m}$. It is evident that a number $m$ produced by such an
algorithm on input string $x$ should be between $0$ and $-\log_2
P(x)$. Let us consider the complementary value $(-\log_2 P(x))-m$
and call it the ``complementary test'' (the consistency requirement
can be easily reformulated for complementary tests).

\medskip

\textbf{Theorem}. \emph{The logarithm of a priori probability
\textup[on the binary tree\textup] $-\log_2 R(x)$ is a
$P$-consistent complementary test for every measure $P$ and has
the usual algorithmic properties.}

\smallskip

Let me remind you that by a priori probability I mean the
universal semicomputable measure introduced in our article with
Zvonkin. [See~\cite{zvonkin-levin}.] It is shown there that it
[its minus logarithm] is numerically close to complexity.

Let us consider a specific computable measure $P$. Compared to
the universal Martin-L\"of test $f$ (specific to a given measure
$P$) our test is not optimal up to an additive constant, but is
asymptotically optimal. Namely, if the universal Martin-L\"of test
finds a deviation $m$, our test finds a deviation at least
$m-2\log_2 m - c$. Therefore, the class of random infinite banry
sequences remains the same.

\smallskip

Now look how nice it fits the philosophy. We say that a
hypothesis ``$x$ appeared randomly according to measure $P$''
can be rejected with certainty $m$ if the measure $P$ is much
less consistent with the appearence of $x$ than a priori
probability (this means simply that $P(x) < R(x)/2^m$. This
gives a law of probability theory that is violated with
probability at most $2^{-m}$. Its violation can be established
effectively since $R$ is [lower] semicomputable [=enumerable from below].
But if this law holds, all other laws of probability theory
[i.e., all Martin-L\"of tests] hold, too. The drawback is that it
gives a bit smaller value of randomness deficiency (only
$m-2\log_2 m -c$ instead of $m$), but this is a price for the
universality (arbitrary probability distribution). The
connection with complexity is provided because $-\log_2 R(x)$
almost coincides with complexity of $x$. Now this connection
does not depend on measure.

It is worth noting that the universal semicomputable measure has
many interesting applications besides the above mentioned. You
know its application to the analysis of randomized algorithms.
Also it is ofter useful in proofs (e.g., in the proof of
J.T.Schwartz' hypothesis regarding the complexity of almost all
trajectories of dynamic systems). Once I used this measure
to construct a definition of intuitionistic validity. All this
show that it is a rather natural quantity.

\medskip

L.

\subsection*{II.}

Dear Andrei Nikolaevich!

I would like to show that plain complexity does not work if we
want to provide an \emph{exact} definition of randomness, even
\emph{for a finite case}. For the uniform distribution on
strings of fixed length $n$ the randomness deficiency is defined
as $n$ minus complexity. For a non-uniform distribution length
is replaced by minus the logarithm of probability.

It turns out that even for a distribution on a finite set the
randomness deficiency could be high on a set of large measure.

\medskip

\textbf{Example}. Let
        $$
P(x)=\begin{cases}
    2^{-(l(x)+100)}, \text{ if } l(x)\le 2^{100};\\
    0, \text{ if } l(x)>2^{100}.
\end{cases}
        $$
Then $|\log_2 P(x)| - K(x)$ exceeds $100$ \emph{for all} strings
$x$.

\smallskip

A similar example can be constructed for strings of some fixed
length (by adding zero prefixes). The violation could be of
logarithmic order.

\smallskip

Let me show you how to sharpen the definition of complexity to
get an exact result (both for finite and infinite sequences).

\textbf{Definitions.}
Let $A$ be a monotone algorithm, i.e., for every $x$ and every $y$ that is
a prefix of $x$, if $A(x)$ is defined, then $A(y)$ is defined too and
$A(y)$ is a prefix of $A(x)$. Let us define
        $$
KM_A(x) = \begin{cases}
        \min\ l(p) \colon x \text{ is a prefix of } A(p);\\
        \infty, \text{ if there is no such } p
          \end{cases}
        $$
The complexity with respect to an optimal algorithm is denoted
by $KM(x)$.

Let $P(x)$ be a computable distribution on the Cantor space
$\Omega$, i.e., $P(x)$ is the measure of the set $\Gamma_x$ of
all infinite extensions of $x$.

\medskip

\textbf{Theorem 1}.
        $$
KM(x) \le |\log_2 P(x)| + O(1);
        $$

\medskip

\textbf{Theorem 2.}
        $$
KM((\omega)_n) = |\log_2 P((\omega)_n)| + O(1)
        $$
\emph{for $P$-almost all $\omega$; here $(\omega)_n$ stands for
$n$-bit prefix of $\omega$. Moreover, the probability that the
randomness deficiency exceeds $m$ for some prefix is bounded by
$2^{-m}$.}

\medskip

\textbf{Theorem 3.}
\emph{The sequences $\omega$ such that}
        $$
KM((\omega)_n) = |\log_2 P((\omega)_n)| + O(1);
        $$
\emph{satisfy all laws of probability theory
\textup(all Martin-L\"of tests\textup).}

\bigskip

Let me use this occasion to tell you the results from my talk in
the laboratory [of statistical methods in Moscow State
University]: why one can omit non-computable tests (i.e., tests
not definable without a strong language).

\smallskip

For this we need do improve the definition of complexity once more.
The plain complexity $K(x)$ has the following property:

\textbf{Remark}. Let $A_i$ be an effectively given sequence of
algorithms such that
        $$
K_{A_{i+1}}(x) \le K_{A_i(x)}
        $$
for all $i$ and $x$. Then there exists an algorithm $A_0$ such
that
        $$
K_{A_0}(x) = 1 + \min_i K_{A_i}(x).
        $$

Unfortunately, it seems that $KM(x)$ does not have this
property. This can be corrected easily. Let $A_i$ be an
effective sequence of monotone algorithms with finite domain
(provided as tables) such that
        $$
KM_{A_{i+1}}(x) \le KM_{A_i(x)}
        $$
for all $i$ and $x$. Let us define then
        $$
\overline{KM}_{A_i}(x) = \min_i KM_{A_i}(x).
        $$
Among all sequences $A_i$ there exists an optimal one, and the
compexity with respect to this optimal sequence is denoted by
$\overline{KM}(x)$. This complexity coincides with the logarithm
of an universal semicomputable semimeasure [=a priori
probability on the binary tree].

\medskip

\textbf{Theorem 4}. \emph{$\overline{KM}(x)$ is a minimal
semicomputable \textup[from above\textup] function that makes
Theorem~$2$ true.}

Therefore no further improvements of $\overline{KM}$ are
possible.

Now consider the language [=set] of all functions computable
with a fixed noncomputable sequence [oracle] $\alpha$. Assume
that $\alpha$ is complicated enough, so this set contains the
characteristic function of a universal enumerable set
[$\mathbf{0}'$].

We can define then a relativized [``языковую'' in the Russian
original] complexity $\overline{KM}_\alpha(x)$ replacing
algorithms by algorithms with oracle $\alpha$, i.e., functions
from this language.

\textbf{Definition}. A sequence $\omega$ is called \emph{normal}
if
        $$
\overline{KM}((\omega)_n) = \overline{KM}_\alpha((\omega)_n)+O(1).
        $$
For a finite sequence $\omega_n$ we define the ``normality
deficiency'' as
        $$
\overline{KM}(\omega_n) - \overline{KM}_\alpha(\omega_n).
        $$

\medskip

\textbf{Theorem 5.} \emph{A sequence obtained by an algorithm
from a normal sequence is normal itself.}

\textbf{Theorem 6.} \emph{Let $P$ be a probability distribution
that is defined \textup(in a natural encoding\textup) by a
normal sequence. Then $P$-almost every sequence is normal.}

This theorem exhibits a law of probability theory that says that
a random process cannot produce a non-normal sequence unless
the probability distribution itself is not normal. This is a
much more general law than standard laws of probability theory
since it does not depend on the distribution. Moreover, Theorem~5
shows that this law is not restricted to probability theory and
can be considered as a univeral law of nature:

\smallskip

\textbf{Thesis}. Every sequence that appears in reality (finite
or infinite) has normality deficiency that does not exceed the
complexity of the description (in a natural language) of how it
is physically produced, or its location etc.

It turns out that this normality law (that can be regarded as
not confined in probability theory) and the law corresponding to
the universal computable test together imply any law of
probability theory (not necessary computable) that can be
described in the language. Namely,the following result holds:

\textbf{Theorem 7}. \emph{Let $P$ be a computable probability
distribution. If a sequence $\omega$ is normal and passes the
universal computable $P$-test, then $\omega$ passes any test
defined in our language \textup(i.e., every test computable with
oracle~$\alpha$\textup).}

{\small

Note that for every set of measure $0$ there exists a test
(not necessary computable) that rejects all its elements.

}

Let us give one more iunteresting result that shows that all
normal sequences have similar structure.

\textbf{Theorem 8}. \emph{Every normal sequence can be obtained
by an algorithm from a sequence that is random with respect to
the uniform distribution.}

\subsection*{III.}

(This letter has no salutation. Levin recalls that he often
gave notes like this to Kolmogorov, who rarely had much
time to hear lengthy explanations and preferred something
written in any case.)

We use a sequence $\alpha$ that provides a ``dense'' coding of a
universal [recursively] enumerable set. For example, let $\alpha$ be
the binary representation of [here the text ``the sum of the a
priori probabilities of all natural numbers'' is crossed out and
replaced by the following:] the real number
        $$
\sum_{p\in A} \frac{1}{p\cdot \log^2 p}
        $$
where $A$ is the domain of the optimal algorithm.

A binary string $p$ is a ``good'' code for $x$ if the optimal
algorithm converts the pair $(p, K(x))$ into a list of strings
that contains $x$ and the logarithm of the cardinality of this list
does not exceed $K(x)+3\log K(x)-l(p)$. (The existence of such a
code means that $x$ is ``random'' when $n\ge l(p)$.)

We say that a binary string $p$ is a canonical code for $x$ if
every prefix of $p$ either is a ``good'' code for $x$ or is a
prefix of $\alpha$, and $l(p)=K(x)+2\log K(x)$.

\textbf{Theorem 1}. \emph{Every $x$ \textup(with finitely many
exceptions\textup) has a canonical code $p$, and $p$ and $x$ can
be effectively transformed into each other if $K(x)$ is given.}

Therefore, the ``non-randomness'' in $x$ can appear only due to
some very special information (a prefix of $\alpha$) contained in
$x$. I cannot imagine how such an $x$ can be observed in
(extracted from) the real world since $\alpha$ is not computrable.
And the task ``to study the prefixes of a specific sequence
$\alpha$'' seems to be very special.

\end{document}